\documentclass
{amsart}

\usepackage[latin1]{inputenc}
\usepackage{subfigure}
\usepackage[english]{babel}
\usepackage{amsmath}
\usepackage{amsfonts}
\usepackage{hyperref}
\usepackage{url}
\usepackage{amssymb}
\usepackage{fullpage}
\usepackage{shuffle}
\usepackage{cases}
\usepackage{tikz}
\usepackage{stmaryrd} 
\usepackage{todonotes}
\usepackage{mathtools}
\usepackage{mathrsfs}
\usepackage{calligra}
\newtheorem{example}{Example}[section]

\newtheorem{theorem}{Theorem}[section]

\newtheorem{lemma}{Lemma}[section]

\def\A{\mathcal{A}}

\def\g{{\mathfrak g}}

\def\N{{\mathbb N}}

\def\U{\mathcal{U}}

\def\al{\alpha}
\def\be{\beta}

\def\scal#1#2{\langle #1 | #2 \rangle}

\def\ncs#1#2{#1\langle\langle #2\rangle\rangle}

\def\ncp#1#2{#1\langle #2\rangle}

\def\Endo{{\mathcal E}\!{\rm nd}}
\def\Lyn#1{{\mathfrak L}{\rm yn}(#1)}
\def\Prim#1{\mathcal{P}{\rm rim}(#1)}

\gdef\stuffle{\;%
  \setlength{\unitlength}{0.0125cm}%
  \begin{picture}(20,10)(220,580) 
  \thinlines 
  \put(220,592){\line( 0,-1){ 10}} 
  \put(220,582){\line( 1, 0){ 20}} 
  \put(240,582){\line( 0, 1){ 10}} 
  \put(230,592){\line( 0,-1){ 10}} 
  \put(225,587){\line( 1, 0){ 10}} 
  \end{picture}\; 
}

\def\adots{\mathinner{\mkern2mu\raise1pt\hbox{.}
\mkern3mu\raise4pt\hbox{.}\mkern1mu\raise7pt\hbox{.}}}

\def\up#1{\raise 1ex\hbox{\footnotesize#1}}

\def\mref#1{(\ref{#1})}

\def\adots{\mathinner{\mkern2mu\raise1pt\hbox{.}
\mkern3mu\raise4pt\hbox{.}\mkern1mu\raise7pt\hbox{.}}}

\def\up#1{\raise 1ex\hbox{\footnotesize#1}}

\def\mref#1{(\ref{#1})}

\usepackage[round]{natbib}

\begin{document}
\author
{M. Deneufch\^atel
}

\title
{Combinatorial Semigroup Bialgebras}

\address[Matthieu Deneufch\^atel]{LIPN - UMR 7030 du CNRS, 99, avenue Jean-Baptiste Cl\'ement, Universit\'{e} Paris 13, 93430 Villetaneuse, France}
\maketitle
\begin{abstract}
This paper is devoted to the presentation of combinatorial bialgebras whose coproduct is defined with the help of a commutative semigroup. We consider this setting in order to give a general framework which admits as special cases the shuffle and stuffle algebras. We also consider the problem of constructing pairs of bases in duality and present Sch\"utzenberger's factorisations as an application.

\noindent {\sc R\tiny{\'ESUM\'E.}}
Cet article est consacr\'e \`a la pr\'esentation de big\`ebres combinatoires dont le coproduit est d\'efini \`a l'aide d'un semi-groupe commutatif. Nous consid\'erons ces structures dans le but de donner un cadre g\'en\'eral dont les alg\`ebres de m\'elange et de quasi-m\'elange sont des cas particuliers. Nous abordons le probl\`eme de la construction de paires de bases en dualit\'e dans ces alg\`ebres; les factorisations de Sch\"utzenberger en sont une application.
\end{abstract}
\section{Introduction}
\label{sec:in}
Recently, stuffle-type products appeared in Combinatorial Physics (\cite{novelli:hal-00632273}). In this paper, we discuss a common generalization of the shuffle and stuffle products. Given an alphabet $X$, the shuffle product, denoted by $\shuffle$, is defined as a law on the algebra $\ncp{k}{X}$ of noncommutative polynomials, bilinear and computed on the words by:
\begin{equation}
\label{shuffle}
        \begin{aligned}
                1 \shuffle w = w \shuffle 1 & = w; \\
                (a u) \shuffle (b v) & = a \Big( u \shuffle (bv) \Big) + b \Big( (au) \shuffle v \Big) 
        \end{aligned}
\end{equation}
for all $u, \, v, \, w \in X^*$ and $a,b, \in X$.\\
The shuffle product appears in several contexts (Lie projectors (\cite{Reutenauer93}), Young tableaux and Littlewood-Richardson rule (\cite{Lothaire02}), noncommutative symmetric functions (\cite{NCSF3}) to cite only a few). One of these contexts is the realm of iterated integrals as shown by Chen's lemma (\cite{Minh2000273,JGJY1}). 

The applications of these products to number theory are well-established. The iterated integrations of the differential forms $\displaystyle \frac{dz}{1-z}$ and $\displaystyle \frac{dz}{z}$ give birth to a family of functions known as \textit{polylogarithms}. The \textit{convergent polylogarithms} are obtained when the first integral is done \textit{w.r.t.} $\displaystyle \frac{dz}{1-z}$. The Taylor expansion of these functions is interesting because its coefficients are series of the form
\begin{equation}
\sum_{n_1> n_2 > \ldots > n_r>0} \, \, \frac{z^{n_1}}{n_1^{k_1}n_2^{k_2}\ldots n_r^{k_r}}
\end{equation}
whose multiplicative structure is given by a perturbation of the shuffle product\footnote{For a discussion about the nuance between deformation and perturbation of the shuffle product, see \cite{SLC62}.} called \textit{stuffle product} (denoted below by $\stuffle$). The definition of this product requires an alphabet \[Y = \left\{ y_n \right\}_{n \geq 1};\] 
for all $y_i, \, y_j \in Y$ and for all $u , \, v \in Y^*$,
\begin{equation}
  \left\lbrace 
\begin{array}{cc}
u \stuffle 1 = 1 \stuffle u = u; \vspace{0.20cm} \\
y_i u \stuffle y_j v = y_i \Big( u \stuffle y_j v \Big) + y_j \Big( y_i u \stuffle v \Big) + y_{i+j} \Big( u \stuffle v \Big).
\end{array}
\right.
\end{equation}
The stuffle product together with the deconcatenation coproduct endows $\ncp{k}{Y}$ with the structure of a graded (commutative) bialgebra (\cite{Hoffman}). This product defines a law of algebra (coassociative with counit) whose value on the letters is obtained as a pertubation of the value of the coproduct dual to the shuffle product:
\begin{equation}
\label{Deltastuffle}
    \Delta_{\stuffle}(y_r) = y_r \otimes 1 + 1 \otimes y_r + \sum_{\genfrac{}{}{0pt}{}{p+q = r}{p,q \geq 1}} y_p \otimes y_q.
\end{equation}

The perturbation term (the sum over $p,q$) of this coproduct is very similar to that which appears in some coproducts used in combinatorial physics (\cite{SLC62}). In this context, the sum is over all decompositions of an element of a semigroup: if $S$ is a semigroup with the finite decomposition property \mref{finitedecomposition} (see \ref{Coproductandcounit}), 
\begin{equation}
   \Delta_{S}: \left\{ \begin{array}{ccc} k \left[ S \right] & \rightarrow & k \left[ S \right] \otimes k \left[ S \right] \\ s & \mapsto & \displaystyle \sum_{rt=s} \, y_r \otimes y_t.
\end{array}\right.
\end{equation}

In this paper, we consider this idea in the general framework of a commutative semigroup (with possibly a zero). The commutativity allows us to prove and apply in this general framework Radford's theorem which is useful when one considers the questions of duality between different bases.

The paper is organized as follows. The next section presents some tools which we then use extensively. Section \ref{generalframework} gives the general framework. In section \ref{shufflealgebra}, we consider the special case of the shuffle algebra in order to illustrate the problem of dual bases. Finally, we present an interesting combinatorial application in section \ref{applications}.

Throughout the paper, $k$ denotes a field of characteristic $0$. By algebras, it will be understood $k$-associative algebra with unit.

\section{Tools}
\label{tools}
\subsection{Multiindex notation}
\label{multiindex}
Let $I$ be a set totally ordered by $<$. If $F=(f_i)_{i\in I}$ is a (totally ordered) family in an algebra $\A$ and if $\al\in \N^{(I)}$ is a multiindex, one defines $F^{\al}$ by
\begin{equation}
        f_{i_1}^{\al_{i_1}} f_{i_2}^{\al_{i_2}} \cdots f_{i_k}^{\al_{i_k}}
\end{equation}
where $J$ is any subset $\{i_1,i_2, \dots , i_k\}\ , \ i_1 > i_2 > \cdots > i_k$, of $I$ which contains the support of $\al$ (it is easily shown that the value of $F^\al$ does not depend on the choice of $J\supset {\rm supp}(\al)$).\\
In particular, if $(e_i)_{i \in I}$ denotes the canonical basis (given by $e_i(j)=\delta_{ij}$) of $\N^{(I)}$, one has $F^{e_i} = f_i$. Morevoer, we set
\begin{equation}
(\alpha + \beta )_i = \alpha_i + \beta_i \, \text{ and } \, \alpha ! = \displaystyle \prod_{i \in {\rm supp}(\alpha)} \alpha_i !, \quad \quad \forall \, \alpha, \, \beta \in \N^{(I)}.
\end{equation}
Finally, we say that a family $\left(T_{\alpha}\right)_{\al \in \N^{(I)}}$ of elements of ${\mathcal A}$ is \textit{multiplicative} (for the product $\times$ of $\mathcal A$) if, $\forall \, \al \, \in \N^{(I)}$, \begin{equation} 
\label{huit}
T_\al \times T_\beta = T_{\al + \be}.
\end{equation}
Note that if $\A$ is commutative, every family is commutative and if $\be$ and $\al$ are such that $\forall i , j \in {\rm Supp}(\al) \times {\rm Supp}(\be), \, \, i \leq j$, then formula \mref{huit} holds.
\subsection{Poincar\'e-Birkhoff-Witt Theorem}
Let $\g$ be a $k$-Lie algebra and $B=(b_i)_{i\in I}$ an ordered basis of $\g$. The Poincar\'e-Birkhoff-Witt (henceforth denoted by PBW) theorem allows one to construct a basis of the enveloping algebra $\U(\g)$.

\begin{theorem}
\label{PBWT}
The elements $B^\al$, for $\al\in \N^{(I)}$, form a basis of $\U(\g)$.
\end{theorem}
The basis obtained by this construction is called a \textit{PBW basis} of $\U(\g)$. This theorem plays a key role in the case of the free algebra as shown below.

\section{Semigroups and Bialgebras}
\label{generalframework}
\subsection{Coproduct and counit}
\label{Coproductandcounit}
Let $(S,\cdot)$ be a commutative semigroup. Define an alphabet $Y$ as follows:
\begin{equation}
 Y = 
 \left\{ 
 \begin{array}{cr}
 \left\{ y_s \right\}_{s \in S \backslash \left\{ \omega \right\}} & \text{ in case } \omega \text{ is a zero of } S; \\  
 \left\{ y_s \right\}_{s \in S} & \text{ otherwise.}
 \end{array} 
 \right.
\end{equation}


As usual, we denote by $Y^*$ the free monoid over $Y$ whose identity element is $1_{Y^*}$.
The non commutative polynomials $\ncp{k}{Y}$ over $Y$ form an algebra for the concatenation product; $\ncs{k}{Y}$ is the set of non commutative series over $Y$ and $\ncp{k}{Y \otimes Y}$ and $\ncs{k}{Y \otimes Y}$\footnote{Alternatively, one could consider $Y^* \otimes Y^*$ as a free partially commutative monoid (see (\cite{DK2})).} denote respectively the set of double polynomials and series.

Define an application $\Delta_{S}: \ncp{k}{Y} \rightarrow \ncs{k}{Y \otimes Y}$ as a morphism of algebras given on the letters by
\begin{equation}
\label{Coproduct}
\Delta_{S} (y_s) = y_s \otimes 1_{Y^*} + 1_{Y^*} \otimes y_s + 
\sum_{
s_1 \cdot s_2 = s}
y_{s_1} \otimes y_{s_2}.
\end{equation}
We say that $S$ has the \textit{finite decomposition property} if, $\forall s \in S \backslash \left\{ \omega \right\}$,
\begin{equation}
\label{finitedecomposition}
      \big| \left\{ (s_1 , s_1), \, s_1 \cdot s_2 = s \right\} \big| < \infty.\tag{D}
\end{equation}
If this is the case\footnote{It is, in particular, the case of:
\begin{itemize}
   \item the semigroup $S=X \cup \{\omega\}$ in which all products equal $\omega$; in this case, one gets the coproduct associated with the shuffle product;
   \item $(\N^+,+)$ (here in its multiplicative notation); in this case, stuffle algebra;
   \item more generally, $(\N^{(X)},+)$ for arbitrary $X$ (finite or infinite).
\end{itemize}
}, the image of $\Delta_{S}$ is an element of $\ncp{k}{Y \otimes Y} \simeq \ncp{k}{Y} \otimes \ncp{k}{Y}$. In fact, $\Delta_{S}$ defines a coassociative coproduct on $\ncp{k}{Y}$.

To complete the bialgebra structure on $\ncp{k}{Y}$, we define the morphism $\epsilon_Y$ whose action on the words is
\begin{equation}
\epsilon_{Y} (w) = 
\left\{ 
\begin{aligned} 
1 & \text{ if } w = 1_{Y^*}; \\ 
0 & \text{ otherwise.} 
\end{aligned} 
\right.
\end{equation}
\begin{lemma}
The map $\epsilon_Y$ is a counit as:
\begin{equation}
   \mu \circ ( \epsilon_{Y} 1_{Y^*} \otimes {\rm Id}) \circ \Delta_{S} = \mu \circ ( {\rm Id} \otimes \epsilon_{Y} 1_{Y^*} ) \circ \Delta_{S} = {\rm Id}.
\end{equation}
\end{lemma}
Therefore, $\mathcal{B} = ( \ncp{k}{Y} , {\rm conc}, 1_{Y^*}, \Delta_{S} , \epsilon_{Y})$ is a bialgebra. It is graded if $S$ is graded.

Assume that $S$ be $\N$-graded, \textit{i.e.} that there exists a function $\ell_S : \, S \rightarrow \N$ such that $S_m \cdot S_n \subset S_{m+n}$, $\forall m , n \in \N^*$ (with $S_n = \left\{ s \in S , \, \ell_S(s) = n \right\}$); assume also that $\ell_S^{-1}(0) = \emptyset$. Then $\mathcal B$ is graded as a bialgebra (\cite{Bourbakalg}). Indeed, one can define $|\cdot| : \, {\mathcal B} \rightarrow \N$ by
\begin{equation}
      |y_s| = \ell_S(s); \quad \quad \displaystyle |y_{s_1} \dots y_{s_k}| = \sum_{i=1}^{k} |y_{s_{i}}|
\end{equation}
and check that the product and coproduct are compatible with the homogeneous components defined by \[{\mathcal B}_p = {\rm span} \displaystyle \left( y_{s_1} \dots y_{s_k}, \, |y_{s_1} \dots y_{s_k}| = p \right), \quad p \in \N.\]

\subsection{Bases}
The letters $y_s$ need not be primitive (for example, for $s \geq 2$, they are not primitive for the stuffle bialgebra; see \mref{Deltastuffle}
). In order to initialize the construction to primitive elements, one has to consider projections of the letters on the set $\Prim{\mathcal{B}}$ of primitive elements of ${\mathcal B}$; this is done with the following projector. \\
From now on we assume that $S$ is graded \textit{with finite fibers} which means that $\displaystyle S = \bigsqcup_{n \in \N} S_n$ with
\[
\left\{
\begin{array}{lr}
   S_0 = \emptyset; & \\
|S_k|<\infty, & \quad \forall \, k \in \N;  \\
S_k \cdot S_\ell \subset S_{k + \ell}.
\end{array}
\right.
\]

We denote by ${\mathcal B}_+$ the kernel of $\epsilon_Y$. One has ${\mathcal B} = {\mathcal B}_+ \oplus k 1_{Y^*}$. Let $I_+$ denote the projector on ${\mathcal B}_+$ along $k$. As $S$ is locally finite\footnote{A locally finite semigroup $S$ is such that $\forall s \in S \backslash \left\{ \omega \right\}$, \[ \Big| \bigcup_{k \geq 1} D_k(s) \Big| < \infty \text{ where } D_k(s) = \left\{ s_1 , \dotsc , s_k \in S \text{ such that } s_1 \cdot \dotsc \cdot s_k = s \right\}.\]}
, $I_+$ is (locally) nilpotent. Therefore, it is possible to define $\pi_1$ by:
\begin{equation}
   \pi_1 = \sum_{k \geq 1} \frac{(-1)^{k-1}}{k} I_+^{*k} = \log_* ({\rm Id}_{\mathcal B})
\end{equation}
where $*$ denotes the convolution product of $\Endo({\mathcal B})$. This operator is a projector on $\Prim{\mathcal{B}}$ (\cite{CartierHopf,Patras98higherlie,Patras}).

\begin{example}
As an example, we consider the bialgebra obtained with the previous setting for $S = (\N^+ , + )$ (we recall that it is the \emph{stuffle algebra} $( \ncp{k}{Y} , {\rm conc} , 1_{X^*} , \Delta_{\stuffle} , \epsilon)$ whose coproduct is given by {\rm Eq.} \mref{Deltastuffle}) 
and compute $\pi_{1}(y_j)$ in this case.

Define, for $b \in {\mathcal B}_+$, $\displaystyle \Delta_+(b) = \Delta(b) - b \otimes 1 - 1 \otimes b$. One can easily check that
\begin{equation}
\label{form4}
   \Big( (I_+\otimes I_+)\circ \Delta \Big)(b) = \Big( (I_+\otimes I_+) \circ \Delta \circ I_+ \Big) (b)  =(\Delta_+ \circ I_+) (b)
\end{equation}
This implies that $\Delta_+$ is a coassociative coproduct on ${\mathcal B}_+$  (since $\Delta$ is coassociative). Hence, on ${\mathcal B}_+$ and for all $k\geq 2$,
\begin{equation}
   \Delta_+^{(k-1)} = I_+^{\otimes k} \circ \Delta^{(k-1)}.
\end{equation}

This relation allows us to compute the convolution powers of $I_+(y_j)$:
\begin{equation}
   I_+^{*k}(y_j) = \mu^{(k-1)} \circ (I_+^{\otimes k}) \circ \Delta^{(k-1)}(y_j) = \mu^{(k-1)} \circ \Delta_+^{(k-1)}(y_j)
\end{equation}
and to prove that
\begin{equation}
   \pi_1(y_j) = \sum_{n > 0} \frac{(-1)^{n-1}}{n} \sum_{\genfrac{}{}{0pt}{}{i_1 + \dotsc + i_{n} = j}{i_1 , \dots , i_{n}>0}} y_{i_1} \dots y_{i_{n}}.
\end{equation}
\end{example}
A total order $<$ on the letters $y_s$ being given
, it is possible to define Lyndon words (whose set will be denoted by $\Lyn{Y}$) and the standard factorisation $\sigma (w)$ of each word $w$ ($\sigma(w)$ is a pair of Lyndon words $\ell_1$ and $\ell_2$ such that $w = \ell_1 \ell_2$ and $\ell_2$ is of maximal length among all the factorisations of $w$). \\
We construct a basis of $\ncp{k}{Y}$ as follows\footnote{Note that this basis is indeed the Lyndon basis of the ``new'' letters $y'_s = \pi_1(y_s)$.}:
\begin{equation}
\label{P_S}
   P_S(w) = \left\{
\begin{array}{ccc}
   \pi_1(y_s) & \text{ if } w=y_s; \vspace{0.2cm}\\
   \left[ P_S(\ell_1) , P_S(\ell_2) \right] & \text{ if } w \in \Lyn{Y} \text{ and } \sigma(w)=(\ell_1 , \ell_2); \vspace{0.2cm} \\
   P_S(\ell_1)^{\alpha_1} \dots P_S(\ell_k)^{\alpha_k} & \text{ if } \left\{ \begin{array}{cc} w = \ell_1^{\alpha_1} \dots \ell_k^{\alpha_k} \\ \ell_1 > \dots > \ell_k \end{array} \right. .
\end{array}
\right.
\end{equation}
(we recall that $\left[ P_1 , P_2 \right] = P_1 P_2 - P_2 P_1, \, \forall \, P_1 , \, P_2 \in \ncp{k}{Y}$).
\begin{theorem}
The elements $P_S(\ell)$, $\ell \in \Lyn{Y}$, form a basis of $\Prim{\mathcal{B}}$. 
\end{theorem}
First, note that if $P, \, Q \in \Prim{\mathcal{B}}$, $\left[ P , Q \right] \in \Prim{\mathcal{B}}$. Since $\pi_1 ( y_s ) \in \Prim{\mathcal{B}}$ for all $s \in S\backslash \left\{ \omega \right\}$, the inductive construction implies that $P_S(\ell) \in \Prim{\mathcal{B}}$, $\forall \ell \in \Lyn{Y}$.

Moreover, the definition of $P_S(\ell)$ implies that it is homogeneous and that \[ P_S(\ell) = \ell + \sum_{u > \ell} \scal{P_S(\ell)}{u} u;\] hence the number of linearly independant $P_S(\ell)$ for a given weight equals the dimension of the corresponding homogeneous component and then $\big(P_S(\ell)\big)_{\ell \in \Lyn{Y}}$ is a basis of $\Prim{\ncp{k}{Y}}$.

The Poincar\'e-Birkhoff-Witt theorem \ref{PBWT} ensures that the family $(P_S(w))_{w \in Y^*}$ forms a basis of $\ncp{k}{Y}$.

In fact, the factorisation of words into Lyndon words implies that the elements of the previous family are also indexed by the multiindices which define their factorisation:
\begin{equation}
\label{correspondencemulti}
\left\{ 
\begin{array}{c}    
w = \ell_1^{\alpha_1} \dots \ell_k^{\alpha_k} \\ \ell_1 > \dots > \ell_k 
\end{array} 
\right\} 
\longleftrightarrow 
\alpha = \alpha_1 e_{\ell_1} + \dots + \alpha_k e_{\ell_k}.   
\end{equation}
Hence, we can denote the elements of the previous family by $P_S(\alpha), \, \alpha \in \N^{(\Lyn{Y})}$. We adopt this notation below.

Let $\left( T_S(\alpha) \right)_{\al \in \N^{(\Lyn{Y})}}$ denote the family dual to $\left( P_S(\beta) \right)_{\beta \in \N^{{(\Lyn{Y})}}}$. It is defined by \[ \scal{T_S(\alpha)}{P_S(\beta)} = \delta_{\alpha \, \beta}.\] \textit{A priori}, $T_S(\alpha)$ is a series. But since ${\mathcal B}$ is graded in finite dimensions, it is in fact a polynomial.

The elements $T_S(\alpha)$ have the following interesting property:
\begin{equation}
\label{mult}
   T_S(\alpha) * T_S(\beta) = \frac{(\alpha + \beta)!}{\alpha ! \beta !} T_S (\alpha+\beta).
\end{equation}
It implies that the family $\left( T_S(e_\ell) \right)_{\ell \in \Lyn{Y}}$ is a transcendence basis of $(\ncp{k}{Y},*)$. This means that the map
\begin{equation}
\left\{
\begin{array}{ccc}
(k[\Lyn{Y}],\cdot) & \rightarrow & (\ncp{k}{Y},*) \\
   \ell & \mapsto & T_S(e_\ell)
\end{array}
\right.
\end{equation}
is an isomorphism of algebras.

\subsection{Statement of the problem}
The facts presented above lead to the following statement: the dual basis of a PBW basis of $(\ncp{k}{Y},{\rm conc})$ is a transcendence basis of $(\ncp{k}{Y},*)$. \\
An example of transcendence basis is given by the set of Lyndon words $\Lyn{Y}$ (which is not necessarily the basis obtained from the dualization presented above, as shown by the example of the free algebra; see \ref{PandS}). This explains why we are interested in the inverse problem: what happens when one starts with a transcendence basis of $(\ncp{k}{Y},*)$ and considers its dual family? In particular, does the dual basis satisfy the same kind of ``Poincar\'e-Birkhoff-Witt relation'' \mref{P_S} as the basis $(P_S(\alpha))_{\alpha \in \N^{(\Lyn{Y})}}$?
\vspace{0.5cm}
\section{Study of the duality from Radford to PBW (shuffle algebra)}
\label{shufflealgebra}
In this section, we are interested in the special case of the shuffle algebra (see \mref{shuffle}) which fits into the previous setting with a semigroup $S = X \cup \left\{ \omega \right\}$ in which all products equal $\omega$. We start with a presentation of several classical results. Then we illustrate the role of Lyndon words as a transcendence basis. Here, $X = \left\{ x_1 > \dots > x_n > \dots \right\}$ is an alphabet and $\Lyn{X}$ denotes the set of Lyndon words over $X$.

\subsection{Bases $(P_w)_w$ and $(S_w)_w$}
\label{PandS}
A basis of the free Lie algebra is given by the standard bracketings of Lyndon words: for $\ell \in \Lyn{X}$,
\begin{equation}
   P_\ell = \left\{
\begin{array}{ccc}
   \ell & \text{ if } & |\ell|=1; \\
\left[P_{\ell_1} ,P_{\ell_2} \right] & \text{ if } & \ell \in \Lyn{X} \text{ and } (\ell_1 , \ell_2) = \sigma(\ell)
\end{array}
\right.
\end{equation}
(we recall that $\sigma (w)$ denotes the standard factorisation of $w$).

If $w = \ell_{i_1}^{\alpha_1} \dots \ell_{i_n}^{\alpha_n}$ with $\ell_1 > \dots > \ell_n$ belongs to $X^*$, let $P_w  = P_{\ell_1}^{\alpha_1} \dots P_{\ell_n}^{\alpha_n}$. The Poincar\'e-Birkhoff-Witt theorem implies that the family $\left( P_w \right)_{w \in X^*}$ is a basis of the free algebra $\ncp{k}{X}$.

Let $(S_w)_w$ be the family of linear forms in $\left( \ncp{k}{X} \right)^* \simeq \ncs{k}{X}$ dual to $(P_w)_w$. It is defined by the relation
\begin{equation}
   \scal{S_w}{P_u} = \delta_{ w \, u}, \quad \forall \, u \in X^*.
\end{equation}
It is possible to prove that $S_w$ is in fact a polynomial given by
\begin{equation}
   S_w = 
\left\{
\begin{array}{ccc}
   w & \text{ if } & |w|=1; \\
   x S_u & \text{ if } & w=xu \text{ and } w \in \Lyn{X}; \\
   \displaystyle \frac{S_{\ell_{i_1}}^{\, \shuffle \, \alpha_1} \shuffle \dots \shuffle S_{\ell_{i_k}}^{\, \shuffle \, \alpha_k}}{\alpha_1 ! \dots \alpha_k !} & \text{ if } & w = \ell_{i_1}^{\alpha_1} \dots \ell_{i_k}^{\alpha_k}\text{ with } \ell_1 > \dots > \ell_k.
\end{array}
\right.
\end{equation}

Note that a similar construction remains possible in the partially commutative case (\cite{DK2}) when one considers an alphabet $X$ with commutations $\theta \subset X \times X$. 

The correspondance \mref{correspondencemulti} implies that we can denote $P_w$ and $S_w$ by the multiindex corresponding to the Lyndon factorisation of $w$. With this notation, one has
\begin{equation}
   S_\alpha \shuffle S_\beta = \frac{(\alpha + \beta)!}{\alpha ! \beta !} S_{\alpha + \beta}, \quad \forall \alpha, \, \beta \in \N^{(I)}.
\end{equation}
Therefore, this family is multiplicative up to a constant.

\subsection[Lyndon words and their dual basis]{Lyndon words as a transcendence basis and their dual basis}
\label{SprimeBprime}
In this section, we take the problem the other way round: we start with a transcendence basis and consider the dual basis. The results of this section have been partially presented, without proof, in (\cite{ISSAC2012}).
\subsubsection{Construction}
\label{Sprimedetail}
A theorem of Radford (\cite{Radford}) ensures that the Lyndon words $\Lyn{X}$ form a transcendence basis for the algebra $(\ncp{k}{X},\shuffle)$. This means that the products $\Lyn{X}^{\shuffle \alpha}$ for $\alpha \in \N^{(I)}$ are a linear basis of $(\ncp{k}{X},\shuffle)$. \\
This explains why we consider the family:
\begin{equation}
\label{Sbasisprime}
S^{'}_w = 
\left\{   
\begin{array}{ccc}
   \ell & \text{ if } & \ell \in \Lyn{X};\\
\displaystyle \frac{S_{\ell_{i_1}}^{'\, \shuffle \, \alpha_1} \shuffle \, \dots \shuffle \, S_{\ell_{i_k}}^{'\, \shuffle \, \alpha_k}}{\alpha_1 ! \dots \alpha_k !} & \text{ if } & w = \ell_{1}^{\alpha_1} \dots \ell_{k}^{\alpha_k}, \, \displaystyle
\left\{
\begin{array}{c}
\ell_1 > \dots > \ell_k \\
\ell_1,\dots,\ell_k \in \Lyn{X}
\end{array}
\right.
\end{array}
\right.
\end{equation}
which is precisely $\big(\Lyn{X}^{\shuffle \al}\big)_{\al \in \N^{(\Lyn{X})}}$. 

The theorem 6.1 in (\cite{Reutenauer93}) states that $S^{'}_w$ is lower triangular:
\begin{equation}
   S^{'}_w = w + \sum_{u < w} \al_u u
\end{equation}
(the coefficients $\al_u$ being integers).

Thus, as this construction is finely homogeneous, it is possible to construct by duality another family of elements of $\ncp{k}{X}$ which will be denoted by $B^{'}_{w}$. We give below several properties of these elements.

\subsubsection{Characterisations of $B^{'}_\ell$}
The following theorems characterise the elements \[B^{'}_\ell, \quad \ell \in \Lyn{X}.\]
\begin{theorem}
\label{theorem1}
 Let $P$ belong to $\ncp{k}{X}$ and $\ell \in \Lyn{X}$. Then
\begin{equation}
P = B^{'}_\ell \Longleftrightarrow \left\{ 
\begin{array}{c}
                               P = \ell + \displaystyle \sum_{\ell < u} \scal{P}{u} u;\\
                               P \text{ is primitive;} \\
                               \forall \ell_1 \in \Lyn{X}, \, \scal{P}{\ell_1} = \delta_{\ell \, \ell_1}.
\end{array}
\right.
\end{equation}
\end{theorem}

\begin{theorem}
\label{reformulation}
Let $P$ belong to $\ncp{k}{X}$ and $\ell \in \Lyn{X}$. Then
\begin{equation}
 P = B^{'}_\ell \Longleftrightarrow 
\left\{ 
\begin{array}{c}
P \text{ is primitive;} \\
|{\rm supp}(P) \cap \Lyn{X}| = 1; \\
\scal{P}{\ell} = 1.
\end{array}
\right.
\end{equation}
\end{theorem}
If $w \in X^*$ is a word whose Lyndon factorisation is $w = \ell_1^{\alpha_1} \dots \ell_n^{\alpha_n}$, $\ell_1 > \dots > \ell_n$, $N_w = \displaystyle \sum_{i=1}^{n} \alpha_i$.
\begin{lemma}
\label{criterionPBW}
The following properties which characterize the elements $B^{'}_w$ are equivalent:
\begin{itemize}
   \item[i)] $\forall \, w \in X^*$ such that $N_w \geq 2$ (\textit{i.e.} such that $w$ is a product of at least two Lyndon words), ${\rm supp}({B^{'}}^w) \cap \Lyn{X} = \emptyset$;
\vspace{0.2cm}
   \item[ii)] $\displaystyle \left( {B^{'}}^w \right)_{w \in X^{*}} = \left( B^{'}_{w} \right)_{w \in X^*}$;
\vspace{0.2cm}
   \item[iii)] $\scal{{B^{'}}^w}{S^{'}_u} = \delta_{w \, u}$ for all $w,u \in X^*$.
\end{itemize}
\end{lemma}

\subsubsection{Recursive construction}
A method similar to Gram-Schmidt algorithm allows us to construct recursively the elements $B^{'}_\ell$ for $\ell \in \Lyn{X}$ from the $P_{\ell'}$ for words $\ell'$ of the same multidegree as $\ell$. In fact, this method gives a way to construct a family in duality with a given linearly independant family. Here, the construction eliminates recursively the Lyndon words that are in the support of $P_\ell$ except for $\ell$ without adding any other Lyndon word and without changing the coefficient of $\ell$.

Let $\alpha \in \N^{(X)}$ be a multiindice and ${\mathcal L}_\alpha = \left\{ \ell_1 < \dots < \ell_m \right\}$ the set of Lyndon words whose multidegree is $\alpha$. 
\begin{lemma}
The elements $B^{'}_{\ell_k}$ for $\ell \in {\mathcal L}_\alpha$ are given by:
\begin{equation}
 \begin{aligned}
B^{'}_{\ell_m} & = P_{\ell_m}; \\
B^{'}_{\ell_{m-1}} & = P_{\ell_{m-1}} - \scal{P_{\ell_{m-1}}}{\ell_m} B^{'}_{\ell_m}; \\
B^{'}_{\ell_{m-2}} & = P_{\ell_{m-2}} - \scal{P_{\ell_{m-2}}}{\ell_{m-1}} B^{'}_{\ell_{m-1}}  - \scal{P_{\ell_{m-2}}}{\ell_{m}} B^{'}_{\ell_{m}}; \\
\vdots & \\
B^{'}_{\ell_{m-k}} & = P_{\ell_{m-k}} - \sum_{j=1}^{k} \scal{P_{\ell_{m-k}}}{\ell_{m-k+j}} B^{'}_{\ell_{m-k+j}}; \\
\vdots
 \end{aligned}
\end{equation}
\end{lemma}

\subsubsection{Multiplicativity property}
\label{counter}
In this section, we give a counter-example which proves that the basis dual to a transcendence basis is not necessarily a Poincar\'e-Birkhoff-Witt basis. 

Let $X = \left\{ a , b \right\}$ be a two-letter alphabet with $a<b$; $B^{'}$ still denotes the basis dual to the basis obtained from the Lyndon words (defined by \ref{Sbasisprime}). 

It is not difficult to show that
\begin{equation}
 B^{'}_{a^2 b^2} = a^2 b^2 - 2 abab + 2 baba - b^2 a^2.
\end{equation}
Since the Lyndon factorisation of the word $(a^2 b^2)^2$ is $(a^2 b^2 , a^2 b^2)$, one has
\begin{equation}
 {B^{'}}^{ ( a^2 b^2 )^2} = \left( B^{'}_{a^2 b^2} \right)^2.
\end{equation}
Hence, $\scal{ a^2 b^2 abab}{{B^{'}}^{ ( a^2 b^2 )^2}} = -2$. Since $a^2 b^2 abab$ is a Lyndon word, the multiplicativity criterion (\ref{criterion}) implies that the basis $\left( B^{'}_w \right)_{w \in X^*}$ is not multiplicative. Moreover since this counter-example is obtained with a square word, it implies that it remains a counter-example for any order on $X$.

\section{Combinatorial applications: Sch\"utzenberger's factorisations} 
\label{applications}
\subsection{General framework}
Let us present the general framework of Sch\"utzenberger's factorisations. Consider a Lie algebra $\g$ and $B=(b_i)_{i\in I}$ an ordered basis of $\g$. Moreover, let $(B_\al)_{\al\in \N^{(I)}}$ denote the Poincar\'e-Birkhoff-Witt basis obtained with $B$. Then (\cite{MelReu}):
\begin{theorem}
\label{gencase}
If $(S_\al)_{\al\in \N^{(I)}}$ is the basis dual to $(B_\al)_{\al\in \N^{(I)}}$, one has
\begin{equation}
\label{main_fact_thm}
	\sum_{\al\in \N^{(I)}} S_\al\otimes B_\al=\prod^{\rightarrow}_{i\in I} \exp\,(S_{e_i}\otimes B_{e_i}).
\end{equation} 
\end{theorem}
Note that the product used to multiply tensor products is $* \otimes \mu_{\U(\g)}$ where $*$ denotes the convolution product of linear forms belonging to $\U^*(\g)$ (the shuffle in the free case) defined by
\begin{equation}
\scal{S_\al * S_\be}{b_i}  =  \scal{S_\al \otimes S_\be}{\Delta(b_i)}, \, \forall i \in I.
\end{equation}
Here $\Delta$ denotes the coproduct associated to the Hopf algebra structure of $\U(\g)$.
Through the completion $\tilde{\Phi}$ of the mapping 
$$
\Phi: V^* \otimes V \rightarrow \Endo^\text{finite} (V)
$$
which associates to each tensor $f\otimes v\in V^* \otimes V$ the endomorphism $\Phi(f\otimes v): b \mapsto f(b) \cdot v$, one gets a factorisation of the identity:\[ \tilde{\Phi} \left (\sum_{\alpha \in \N^{(I)}} S_\alpha \otimes B_\alpha \right) = {\rm Id_{\U(\g)}}.\]
It is interesting to consider the inverse problem: is it possible to write such a factorisation of the identity when one starts with a multiplicative basis of $\U^*(\g)$ and consider its dual basis?

Let $\left( S_\al \right)_{\al \in \N^{(I)}}$ be a basis of $\U^*(\g)$ and considers its dual family $(B_{\be})_{\be \in \N^{(I)}}$ 
.

On the one hand, the duality of the families $( B_{\al} )_\al$ and $(S_{\al})_\al$ allows us to write a kind of resolution of the identity:
\begin{equation}
 {\rm Id}_{\U(\g)} = \sum_{\al \in N^{(I)}} S_\al \otimes B_{\al}.
\end{equation}

On the other hand, one can consider the product that defines the left-hand side of the factorisation and, assuming that the family $\left( S_\al \right)_{\al \in \N^{(I)}}$ satisfies relation \mref{mult}, one obtains the following derivation
:
\begin{equation}
 \begin{aligned}
  \prod^{\rightarrow}_{i\in I} \exp\,(S_{e_i}\otimes B_{e_i}) & = \sum_{k \geq 0} \sum_{\genfrac{}{}{0pt}{}{i_1 \geq \dots \geq i_k}{\alpha_1 , \dots , \alpha_k}} \frac{(S_{e_{i_1}} \otimes B_{e_{i_1}})^{\alpha_1} \dots (S_{e_{i_k}} \otimes B_{e_{i_k}})^{\alpha_k} }{\alpha_1! \dots \alpha_k!}\\
& = \sum_{k \geq 0} \sum_{\genfrac{}{}{0pt}{}{i_1 \geq \dots \geq i_k}{\alpha_1 , \dots , \alpha_k}} \frac{S_{e_{i_1}}^{\alpha_1} * \dots * S_{e_{i_k}}^{\alpha_k}}{\alpha_1! \dots \alpha_k!} \otimes (B_{e_{i_1}})^{\alpha_1} \dots (B_{e_{i_k}})^{\alpha_k} \\
& = \sum_{\al \in \N^{(I)}} S_\al \otimes \prod_{i \in {\rm supp}(\al)} (B_{e_i})^{\alpha_i} = \sum_{\al \in \N^{(I)}} S_\al \otimes B^\al.
\end{aligned}
\end{equation}

But since we do not know whether the family $\left( B_{\al} \right)_{\al \in \N^{(I)}}$ is of Poincaré-Birkhoff-Witt type, we can not go further: \textit{a priori}, $B^\al = \displaystyle \prod_{i \in {\rm supp}(\al)} (B_{e_i})^{\alpha_i} \neq B_{\al}$ (and we present in section \ref{counter} a counter-example). \\
Moreover, even if one can easily construct a multiplicative family $\left( Q^\al \right)_{\al \in \N^{(I)}}$ (defining the element $Q_{\al}$ by $\displaystyle \prod_{i \in {\rm supp}(\al)} \left( B_{e_i} \right)^{\al_i}$), the family $\left( Q_{\al} \right)_{\al \in \N^{(I)}}$ is \emph{a priori} not dual to $\left( S_{\al} \right)_{\al \in \N^{(I)}}$ and does not lead to such a factorisation.
  
This explains why we want to know in which cases one has the equality between $B^{\al}$ and $B_{\al}$ for every multiindex $\al$.
\subsection{General case: a multiplicativity criterion}
It is possible to give a necessary and sufficient condition for a basis obtained by duality to be multiplicative. Roughly speaking, for the basis to be multiplicative, it must vanish on all the products of more than two elementary elements of its dual basis.
\begin{lemma}
\label{criterion}
Let $\left( S_\alpha \right)_{\alpha \in \N^{(I)}}$ be a basis of ${\mathcal U} ({\mathfrak g})$ satisfying $\scal{S_{e_i}}{1_{\U(\g)}} = 0$, for all $i \in I$, in duality with a basis $\left( B_\alpha \right)_{\alpha \in \N^{(I)}}$ : $\scal{S_\alpha}{B_{\beta}} = \delta_{\alpha \beta}$. Then $B^\beta = \displaystyle \prod_{i \in {\rm supp}(\beta)} B_{e_i} = B_{\beta}$ ($B$ is multiplicative) if and only if 
\begin{equation}
 \forall i \in I, \, \forall \beta \in \N^{(I)}, \, |\beta| \geq 2, \, \, \scal{S_{e_i}}{B^{\beta}} = 0.
\end{equation}
\end{lemma}
\section{Conclusion}
The ideas presented above lead us directly to some unanswered questions. First, one can ask whether it is possible, in the previous criterion, to consider only multiindices $\beta$ with $|\beta| = 2$ instead of $|\beta| \geq 2$. Then one can try to understand the structure of the elements $B^{'}_{w}$. Since they are primitive, they might come from a bracketing process; is it true? If this is the case, what is the good bracketing? Moreover, we know that it is not possible to write Sch\"utzenberger's factorisation with the bases $S^{'}_w$ and $B^{'}_w$. But it is possible to apply the Poincar\'e-Birkhoff-Witt theorem to the family $B^{'}_{\ell}$ for $\ell \in \Lyn{X}$ to construct a basis of the free associative algebra; then by duality, one obtains a new basis which has the good behaviour \textit{w.r.t.} Sch\"utzenberger's factorisation. It would be interesting to compare the elements of this new basis with the elements $S^{'}_w$ in order to get a deeper understanding of the properties that need to be 
satisfied in order to write the factorisation in its \textit{resolution of the identity} form. Note that we used numerical computations to understand the properties of the objects presented in this paper; because of the behaviour of the shuffle product, the interesting examples arise for words of length $\geq 8$ (see \ref{counter} for example), which leads to long computations. Therefore, the numerical investigations are not as easy as one could expect for these questions.

The numerical computations have been carried out with the software Sage (\cite{Sage-Combinat}). The bialgebra structures are implemented in the case of the shuffle and stuffle algebras. These programs are presented in two worksheets available at the following addresses:
\begin{itemize}
   \item \url{http://sagenb.org/home/pub/4504/} (shuffle algebra);
   \item \url{http://sagenb.org/home/pub/4519/} (stuffle algebra).
\end{itemize}

\textbf{Acknowledgements}
\label{sec:ack}
The author wants to thank G\'erard H. E. Duchamp for his advice and V. Hoang Ngoc Minh for fruitful discussions.

\bibliographystyle{abbrvnat}
\bibliography{Biblio}
\label{sec:biblio}

\end{document}